\documentclass[letterpaper, 10 pt, conference]{ieeeconf}  % Comment this line out if you need a4paper
\newcommand{\fref}[1]{(\ref{eq:#1})}
\usepackage{adjustbox}
\usepackage{array}
\usepackage{booktabs}
\usepackage{multirow}
\usepackage{hhline}

\newcolumntype{R}[2]{%
    >{\adjustbox{angle=#1,lap=\width-(#2)}\bgroup}%
    l%
    <{\egroup}%
}
\IEEEoverridecommandlockouts                              % This command is only needed if 
\usepackage{graphicx}
\usepackage{amsmath}
\usepackage{color}

\title{\LARGE \bf
Quasi-dynamic Load and Battery Sizing and Scheduling for Stand-Alone Solar System Using Mixed-integer Linear Programming
}

\author{Abdulelah H. Habib, Vahid R. Disfani, Jan Kleissl and Raymond A. de Callafon
\thanks{Abdulelah Habib, Vahid R Disfani, Jan Kleissl and Raymond A. de Callafon are with the Department of Mechanical and Aerospace Engineering, University of California, San Diego
       {\tt\small {ahhabib,disfani,jkleissl,callafon}}
       @ucsd.edu}}%\\$^*$ Corresponding author}%
\begin{document}

\maketitle
\thispagestyle{empty}
\pagestyle{empty}

%%%%%%%%%%%%%%%%%%%%%%%%%%%%%%%%%%%%%%%%%%%%%%%%%%%%%%%%%%%%%%%%%%%%%%%%%%%%%%%%
\begin{abstract}
Considering the intermittency of renewable energy systems, a sizing and scheduling model is proposed for a finite number of static electric loads. The model objective is to maximize solar energy utilization with and without storage. For the application of optimal load size selection, the energy production of a solar photovoltaic is assumed to be consumed by a finite number of discrete loads in an off-grid system using mixed-integer linear programming. Additional constraints are battery charge and discharge limitations and minimum uptime and downtime for each unit. For a certain solar power profile the model outputs optimal unit size as well as the optimal scheduling for both units and battery charge and discharge (if applicable). The impact of different solar power profiles and minimum up and down time constraints on the optimal unit and battery sizes are studied. The battery size required to achieve full solar energy utilization decreases with the number of units and with increased flexibility of the units (shorter on and off-time).
A novel formulation is introduced to model quasi-dynamic units that gradually start and stop and the quasi-dynamic units increase solar energy utilization. The model can also be applied to search for the optimal number of units for a given cost function.

\end{abstract}

%%%%%%%%%%%%%%%%%%%%%%%%%%%%%%%%%%%%%%%%%%%%%%%%%%%%%%%%%%%%%%%%%%%%%%%%%%%%%%%%
\section{INTRODUCTION}
\def\paraph{\textcolor{blue}}
Standalone solar energy systems are increasingly deployed in rural and off grid areas, especially to provide basic societal needs, for example, water treatment, pumping, and cooking or heating \cite{iea2012renewable}. The main obstacles in optimal sizing and scheduling problems of electric power systems are the variability and intermittency of renewable energy generation. The most challenging scenario is the standalone or islanded mode where high penetration of variable renewable power sources such as wind and solar causes power variability that is large enough to affect electric power quality and efficiency
\cite{nrelvar} and \cite{pvvoltage}%\cite{nrelvar,pvvoltage} and \cite{windsolarchalng}. 
Solar energy production depends foremost on the solar resource availability, which suffers from high variability over a broad range of time scales.

Sizing of power generators for standalone application was discussed in many papers such as \cite{reviewstandalonePV} and \cite{standalonesizing1}, %\cite{reviewstandalonePV,standalonesizing} and \cite{standalonesizing1},
where solar energy was sized to meet the load. In this paper, we consider the opposite case where both the generators, unit load size, and number of units are designed to optimize system efficiency or minimize energy loss. Such a tool is helpful if a complete microgrid (generators and loads) is designed from scratch, such as for a desalination plant without local grid power supply, or if additional loads are connected to an existing off-grid solar system.

Load scheduling plays an important role in optimizing efficiency as well. It has been applied in many ares, for example thermal loads and domestic appliances \cite{c3}. Optimal scheduling a wind farm with a storage system constrained by states of charge of the battery was considered in \cite{c9}. Game theory and customers effect on the grid and EVs was studied in \cite{c2}. %Another example combined domestic appliance scheduling, EV, and wind variability \cite{c4}. Also HVAC system has been included in the scheduling problem \cite{c5}. 
Load scheduling is also used in water network system where the number of pumps are scheduled to meet the water demand and optimize the cost \cite{schedulingpumps}. Forecasts for solar generation and uncontrollable loads are required inputs for the scheduling problem yet solar forecast research is still ongoing and errors can be substantial especially on short time scales. In our previous work we defined number and size of units and solved the scheduling problem with the solar forecast. We also proposed ideas to overcome forecast error in standalone cases \cite{habibsyscon}.

Energy storage systems can solve the variability and intermittency problem \cite{c1} and balance forecast errors, but energy storage will add cost and complexity to the standalone system. An alternative way is load sizing, where loads on the standalone system are adjusted to accommodate power variability to consume as much as possible of the available solar energy thereby reducing energy losses. While the need for energy storage can be substantially reduced through scheduling, the addition of a battery can be cost effective for standalone applications be it wind \cite{windbattery} or solar PV \cite{solarbattery}. The sizing and scheduling of such a battery will also be optimized in this work.

 This paper is organized as follows, some background and problem explanation are discussed in Section~\ref{sec:background}. In Section~\ref{sec:problemformulation} the problem formulation is given explaining the optimization techniques, objective functions and the constraints. Also the two principal optimization approaches are presented in this section. In Section~\ref{sec:results} different scenarios are discussed and compared for different solar generation pattern, e.g., clear, cloudy and partially cloudy days as well as different minimum up and downtime constraints both with and without storage.

\section{Background }
\label{sec:background}
Planning algorithms for solar standalone applications are needed to overcome solar radiation variability. %and hard to predict over years, although the seasonal effect is easier.
The standalone PV sizing problem inputs are solar PV and the load. Energy storage is an optional component. PV power performance models are well understood and PV output timeseries can be generated globally and for different days of the year given existing solar resource databases and models. For brevity we will deal with sample daily solar power profiles from a PV system at UC San Diego and normalized it to one for generality. In practice, sizing decisions should be based on several years of solar resources data to capture the annual variability and possible even inter-annual variability. While typical meteorological year (TMY) or typical solar year are often used for this purpose, large solar system developers increasingly rely on multidecadal modeled power production based on site adaption of long-term satellite records with short-term local measurements for their financial calculations \cite{longtermpred}. Such long-term data would be preferable in practice although interannual variability of solar energy generation is small. For example, \cite{GermanSolar} specifies that the interannual variability of GHI for 7-10 years of measurement at Potsdam, Germany and Eugene, USA is about 5\%.  

The sizing of load is determined to maximize the solar utilization factor that will be referred to in this paper as efficiency (Eff). We proposed a unit sizing design for each a clear and a cloudy day. The approach can be extended to yearly data which allows characterizing most of the important seasonal and diurnal variability.  Considering longer timeseries merely adds computational cost to the implementation of our proposed model. If computation power is limited, the sizing could be based on these two or a few more characteristic days with a weighing factor based on the probability of occurrence of a daily pattern is for a selected location. 

In our previous work we targeted the sizing of standalone PV reverse osmosis units (PVRO) by searching among %\jan{rephrase 'searching among' to better characterize the method}
different unit numbers and sizes. A financial model allowed the optimization to achieve the lowest water cost \cite{habibasme}. Then we developed a mathematical data driven optimization to optimally size the PVRO units in \cite{habibRE}. In \cite{habibsyscon,habibACC}, we developed a model predictive load scheduling to optimally schedule units. In this paper we describe a method that can consider all of these objectives at the same time; we provide a model to optimally size, choose the number of units, and schedule. The addition of energy storage is also considered. The algorithm can also be applied to problems where units sizes and numbers are known and only the scheduling is of interest. %\jan The preceding paragraph would fit better in section I. Also add a sentence about the methdological novelty described in III.A (dynamics of load switching).
\begin {figure}[ht]
\graphicspath{ Plots/ }
\includegraphics[width=.9\columnwidth]{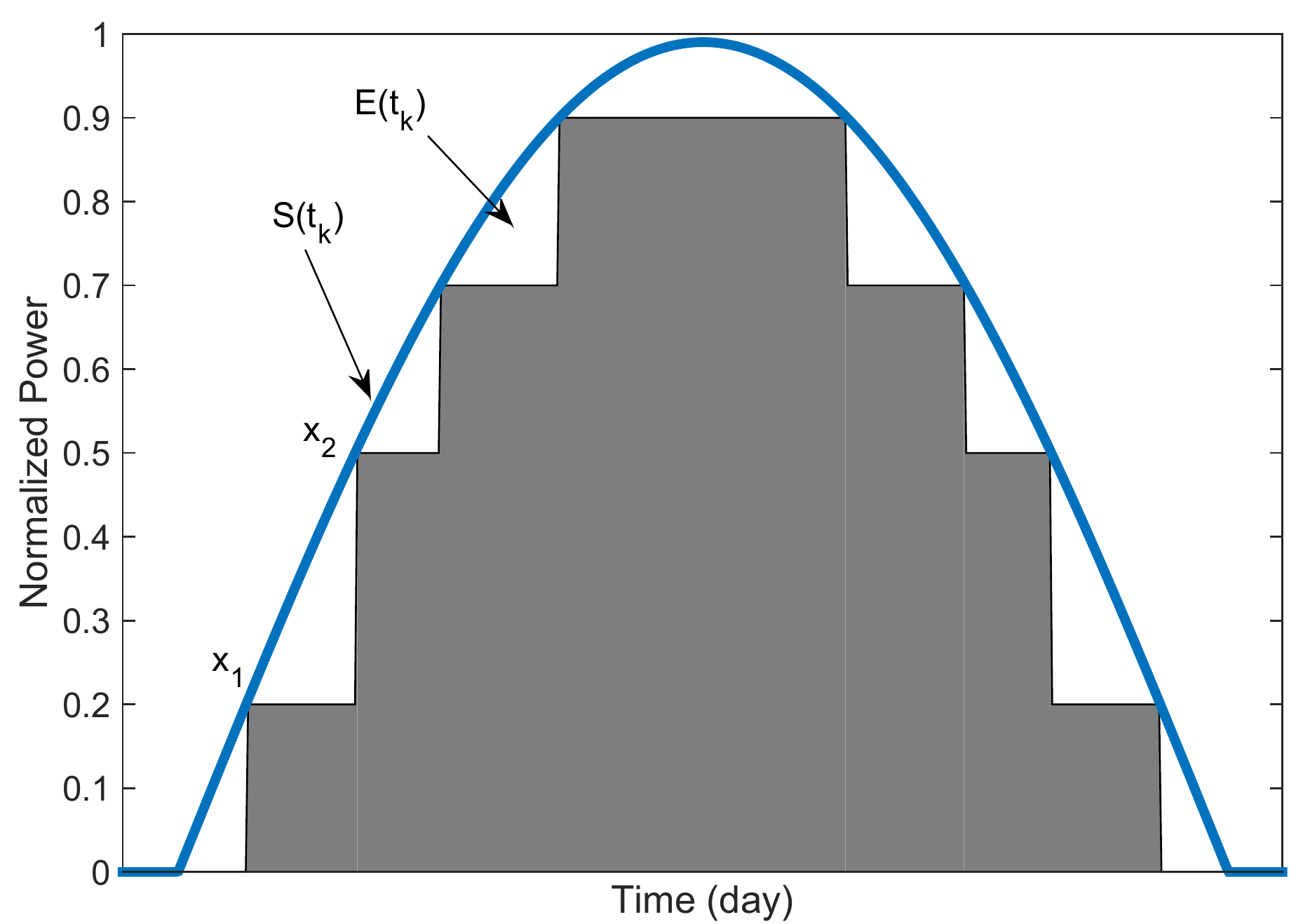}
\centering
\caption{Problem illustration for one clear and symmetric day. Loads $x^i$ are scheduled hourly to follow the increase or decrease in solar power generation $S(t_k)$.}
%\jan specify the x used in this example. If there is only x1 and x2 I do not understand you can obtaint the largest output at 0.9.
\label{figure:min}\vspace{-10pt}
\end {figure}

\section{Problem Formulation }
\label{sec:problemformulation}
The optimization problem tackled in this paper is generally presented as
%Denoting the binary variable $u_i^{t_k}$ as the on/off status of the unit $i$ at time step $t_k$, the optimization problem \fref{opt} is generally presented as below,
% \begin{equation}
% \begin{aligned}
% & \underset{u_i, P_{L_i}, P_b}{\text{minimize}}
% & & \sum\limits_{\substack{t=1}}^{T} [ P_{PV}(t)-\sum\limits_{\substack{i=1}}^{n} y_i(t)]\\
% \end{aligned}
% \end{equation}
\begin{equation}
\begin{aligned}
 \min_\mathbf{U,\text{X},{\bar{\rm P}}_b,{\rm P_s}} &&\mathbf{f(U,\text{X},{\bar{\rm P}}_b,{\rm P_s}})\\
 \text{s.t.}&&\mathbf{g(U,\text{X},{\bar{\rm P}}_b,{\rm P_s})}\le0\\
&& \mathbf{h({U},\text{X},{\bar{\rm P}}_b,{\rm P_s})}=0\\
&& \mathbf{\phi({U},R,Q,W,V)}\le0\\
&& \mathbf{\psi({U},R,Q,W,V)}=0\\
% && \phi:\mathbf{U}\rightarrow\{0,1\}^{2\times|\mathbf{U}|}\\
% &&& \phi:\mathbf{U}\rightarrow\{0,1\}^{2\times|\mathbf{U}|}\\
% && \mathbf{U}=\frac{{R+Q}}{2}\\
% && {R-Q}\ge0\\
&& \mathbf{R,Q,W,V}\subset\{0,1\}^{|\mathbf{U}|}\\
&& \mathbf{U}\subset\{0,1/2,1\}^{|\mathbf{U}|}\\
\end{aligned}
\label{eq:mainopt}
\end{equation}
where the decision variables $\text{X}=[x_1 ,x_2,\cdots,x_n]^T$ and
% \[
% \mathbf{U}=[uk]=[u_1 ,u_2,\cdots,u_N]^T 
% \]
\begin{equation*}
\begin{array}{rcl}
\mathbf{U} &=& 
{\left [ \begin{array}{cccc} 
u^1_{1} & u^2_{1} & \cdots & u^k_{1} \\
u^1_{2} & u^2_{2} & \cdots & u^k_{2} \\ 
\vdots &  & \ddots& \vdots \\
u^1_n & u^2_n & \cdots & u^k_n
\end{array} \right ] %\otimes \mathbf{1}_{1 \times L}} \in R^{T \times n}\\
%x &=& \left [ \begin{array}{cccc} x^1 & x^2 & \cdots & x^n \end{array} \right ]^T \in R^{n \times 1}
}\end{array}
\label{eq:U1}
\end{equation*}
are load vector and switching matrices, respectively, such that  $u_i^{t_k}$ denotes the portion of the load $i$ which is turned on at time $t_k$, and $x_i$ is the size of unit $i$. The size of the battery connected is a scalar and called P$_b$, while P$_s$ is a column vector representing the battery scheduled discharging power at all $t$. The variable matrices $\mathbf{V}$ and $\mathbf{W}$ respectively denote the start-up and shut-down signals, whereas the variable matrices $\mathbf{R}$ and $\mathbf{Q}$ are defined such that $\mathbf{U=\frac{R+Q}{2}}$. Also, the functions $\mathbf{f}$, $\mathbf{g}$ and $\mathbf{h}$ denote the mathematical formulations of the objective function, inequality constraints, and equality constraints respectively. The functions $\phi$ and $\psi$ also represent the inequality and equality constraints that involve $\mathbf{U}$ and the binary matrices.

\subsection{Objective Function}
As indicated in Figure~\ref{figure:min}, the power mismatch between available solar power and power used by load units at any time $t_k$ can be characterized by
 \begin{center}
 $E(t_k)= S(t_k)-\sum_{i=1}^{n} u^i_kx^i$ \end{center}
for the case without a battery. Considering a battery system that is discharged by $P_s^{t_k}$ at time $t_k$, the definition of power mismatch changes to 
\begin{center}
$E(t_k)= S(t_k)+P_s(t_k)-\sum_{i=1}^{n} u^i(t_k)\cdot x^i$
\end{center}
%where $u_k^i \in [0,1]$ is a $n$ binary numbers reflecting the on/off {switch state} of the individual loads $i=1,\ldots,n$ with their {(to be determined) static load size $x^i$}. Defining the vectors
% \begin{equation}
% \begin{array}{rcll}
% u_k &=& \left [ \begin{array}{cccc} u_k^1 & u_k^2 & ... & u_k^n \end{array} \right ],~ & {u_k^i} \in [0,1]\\ x &=& \left [ \begin{array}{cccc} x^1 & x^2 & ... & x^n \end{array} \right ]^T,~ & x^i > 0
% \end{array}
% \label{eq:vars}
% \end{equation}
% the static power mismatch $E(t_k)$ at a particular time $t_k$ can be written with an inner product
% \[
% E(t_k) = S(t_k) - u_k x,
% \]
% of the {time dependent binary switch state vector $u_k$ and the static load size distribution vector $x$}. Static load response optimization can now be written as
% \begin{equation}
% \mbox{arg} \min_{u_k,x} \sum_{k=1}^N E(t_k)^2,~~ E(t_k) = S(t_k) - u_k x 
% \label{eq:opt}
% \end{equation}
% where the variables $u_k$ and $x$ are given in \fref{vars}.

Therefore, the objective function is to minimize the mismatched power and maximize the efficiency (defined as matched energy over total solar energy)
\begin{equation}%\vspace{-3pt}
\begin{aligned}
\mathbf{f(U,\text{X},{\bar{\rm P}}_b,{\rm P_s}})=\mathbf{1}^T\cdot\left(S+{\rm P_s}-diag(\text{X})\cdot U\cdot\mathbf{1}\right)
\end{aligned}
\label{eq:obj1}\vspace{-5pt}
\end{equation}
where $S$ is a column vector denoting the PV power available at all time steps. The objective function in \fref{obj1} is nonlinear due to the bilinear product function, \emph{i.e.} $diag(\text{X})\cdot U$. To remove nonlinearity, a new decision variable matrix $\mathbf{Y}=diag(\text{X})\cdot U$ is defined which denotes the scheduled power. Thus, the objective function becomes linear as
\begin{equation}
\begin{aligned}
\mathbf{f(U,\text{X},{\bar{\rm P}}_b,{\rm P_s}})=\mathbf{1}^T\cdot\left(S+{\rm P_s}-\mathbf{Y}\cdot\mathbf{1}\right)
\end{aligned}
\label{eq:mainopt2}\vspace{-5pt}
\end{equation}

Adding the definition of $\mathbf{Y}$ to the set of constraints guarantees identical solutions for Eqs. \fref{mainopt} and \fref{mainopt2}. The definition of $\mathbf{Y}$ moves nonlinearity from objective function to constraints. The big-M method is the common solution to remove such nonlinearities \cite{bigM}, when $\mathbf{U}$ is a binary matrix. In this paper, a novel application of big-M method is proposed to remove these nonlinearities when the elements of $\mathbf{U}$ belong to the set of $\{0,1/2,1\}$.
\subsection{Constraints}
\subsubsection{Resource Adequacy}
To prevent frequency issues, the maximum total load that the microgrid can supply must be less than the total PV energy available at each time interval:
\begin{equation}
\begin{aligned}
\mathbf{1}^{T} . \mathbf{Y} \le {\rm S} +{\rm P_s},
\end{aligned}
\label{const_power}
\end{equation}
%where $\mathbf{S}=[s_1,s_2,\cdots,s_k]$ is the PV generation vector, and {\rm P$_s$}$=[{p_s}_1,{p_s}_2,\cdots,{p_s}_k]$ is the power of the battery. 
\subsubsection{Definition of Scheduled Power}
In order to apply the new format of big-M method, two auxiliary binary matrices $R$ and $Q$ are defined such that $\mathbf{U=\frac{R+Q}{2}}$ and ${R-Q}\ge0$. These two constraints guarantee that the vector $U=[0,1/2,1]$ is uniquely mapped to the vectors $R=[0,1,1]$ and $Q=[0,0,1]$. 

With the definitions of the matrices $R$ and $Q$, the constraint $\mathbf{Y}=diag(\text{X})\cdot U$ is equivalent to the following set of constraints, 
\begin{equation}
\begin{aligned}
\mathbf{-Y\le0}\\
\mathbf{Y}-\frac{\mathbf{R+Q}}{2}M\mathbf{\le0}\\
\mathbf{Y}-\frac{\mathbf{1\cdot1^T\cdot diag({X})}}{2}-\frac{\mathbf{R+1\cdot1^T-Q}}{2}M\mathbf{\le0}\\
-\mathbf{Y}+\frac{\mathbf{1\cdot1^T\cdot diag({X})}}{2}-\frac{\mathbf{R+1\cdot1^T-Q}}{2}M\mathbf{\le0}\\
\mathbf{Y-1\cdot1^T\cdot diag(X)\le0}\\
\mathbf{-Y+1\cdot1^T\cdot diag(X)}-\frac{\mathbf{1\cdot1^T-R+Q}}{2}M\mathbf{\le0}
\end{aligned}
\label{eq:const_bigM}
\end{equation}where $M$ is a real constant number, {e.g}, $10^6$.

In \eqref{eq:const_bigM} the first two and the last two constraints correspond to the case that $u=0$ and $u=1$, respectively. the third and fourth inequalities guarantee any $\mathbf{y}=x/2$ if the corresponding $u$ is $0.5$, while they are relaxed otherwise. 
 
\subsubsection{Dynamic Model Constraint} One of the novelties in this paper is to consider dynamics of load switching in the start up and shut down processes. The dynamics in this paper is modeled by adding an intermediate step $u=1/2$ for the load units while switching on or off. 

To model these dynamics, several constraints must be considered. First, no immediate transitions between $u=0$ and $u=1$ are allowed. Second, if $u(t)=1/2$, then at the following timestep $u(t+1) = 1-u(t-1)$, which means that $u$ may not stay in the state $u=1/2$ for any two consecutive time steps and no switching is allowed during dynamics. These constraints are summarized as:
\begin{equation}
\begin{aligned}
u_i^{t}-u_i^{t-1}\le1/2 && \forall_{i,t}\\
u_i^{t-1}-u_i^t\le1/2 && \forall_{i,t}\\
r_i^{t+1}\le 1-r_i^{t-1}+(1-r_i^t+q_i^t)M&& \forall_{i,t}\\
r_i^{t+1}\ge 1-r_i^{t-1}-(1-r_i^t+q_i^t)M&& \forall_{i,t}\\
q_i^{t+1}\le 1-q_i^{t-1}+(1-r_i^t+q_i^t)M&& \forall_{i,t}\\
q_i^{t+1}\ge 1-q_i^{t-1}-(1-r_i^t+q_i^t)M&& \forall_{i,t}\\\end{aligned}
\label{const_Y}
\end{equation}

\subsubsection{Minimum Up-time and Minimum Down-time} 
To avoid increased wear and tear to load units and inconvenience to microgrid customers because of frequent start-ups and shut-downs, a set of constraints are defined to guarantee that the unit is switched on (off) for at least $m^+$ ($m^-$) time steps before it is switched off (on). These constraints are called minimum up (down) time and are defined as:
\begin{equation}
\begin{aligned}
r_{i,t_k}-\sum_{h=t_k-m^+_i +2}^{t_k} v_{i,h}\leq 0 && \forall_{m^+_i\le t_k\le T}\\
(1-r_{i,t_k})-\sum_{h=t_k-m^-_i }^{t_k} w_{i,h}\leq 0&& \forall_{m^-_i\le t_k\le T},
\end{aligned}
\label{const_minupdown}
\end{equation}
where the matrix $\mathbf{V}\subset\{0,1\}^\mathbf{|V|}$ and $\mathbf{W}\subset\{0,1\}^\mathbf{|W|}$ are denoted as start-up and shut-down matrices respectively, and their elements are defined as:
\begin{equation}
\begin{aligned}
&v_{i,t_k}-w_{i,t_k}=r_{i,t_k}-r_{i,t_k-1} &&\forall_{1\le i\le N}\forall_{2\le t_k\le T}\\
&v_{i,t_k}+w_{i,t_k}\leq 1 && \forall_{1\le i\le N}\forall_{2\le t_k\le T}\\
&v_{t,1}=w_{i,1}=0&& \forall_{1\le i\le N}
\end{aligned}
\label{const_vw}
\end{equation}

% To differentiate between demand matrix $\mathbf{X}$ and supplied load matrix $\mathbf{Y}={diag}(\mathbf{X}) \times \mathbf{U}$ for more clear representations in the rest of the paper. Now, the nonlinearity is moved to the constraints and it can be resolved by applying big-M method to define the alternative linear constraints as below \cite{bigM}:
% \begin{equation}
% \begin{aligned}
% \mathbf{-Y\le0}\\
% \mathbf{Y-U}M\mathbf{\le0}\\
% \mathbf{Y-1\cdot1^T\cdot diag(X)\le0}\\
% \mathbf{1\cdot1^T\cdot diag(X)+(U-1\cdot1^T)}M\mathbf{-Y\le0}\\
% \end{aligned}
% \label{eq:const_bigM}
% \end{equation}
% y>=0
% y<=u*1e6
% for i=1:nu_loads
%     y(i,:)<=x(i)
%     y(i,:)>=x(i)+(u(i,:)-1)*1e6 %big M method  
% end
% where $M$ is a real constant number, {e.g}, $10^6$.
% The optimization problem considered in this paper can be presented in a general format as below,

% where $i$ and $t_k$ are the index for units and time steps respectively. $S(t_k)$ denotes the solar power timeseries of length \(T\). $x \in R^{n}$ is the vector of load sizes. $y_i(t_k)$ also denotes the committed demand by load $i$ at time $t_k$.

\subsubsection {Battery Constraints}
% \end{itemize}
% \begin{equation}
% \begin{aligned}
% \left\{\begin{matrix}
% &&& P_b(t)=
% \left\{\begin{matrix}
%  0,  & & u_b=0\\
% S(t) -\sum\limits_{\substack{i}} y_i(t),& & u_b=1
%  \end{matrix}\right.\\
% &&& P_b (t) \leq \bar{P}_b, \forall t \\
% &&& P_b (t) \geq -\bar{P}_b, \forall t\\
% &&& E_b(t+1) =E_b(t) + P_b (t) \tau, \forall t\\
% &&& 0 \leq E_b(t) \leq \bar{E_b}, \forall t\\
% &&& E(0)=E_{initial}, E(0)=E(T), \\
% &&& u_b (t)\in \{0,1\},  \forall t \\ \\
% \end{matrix}\right.\\
% \end{aligned}
% \label{eq:MILP1}
% \end{equation}
There are also some constraints associated with the battery such as maximum charge and discharge power, minimum and maximum limits of state of charge (SOC), and initial and final SOC, which are described as:
\begin{equation}
\begin{aligned}
{\rm P_s}\le {\bar{\rm P}}_b\cdot \mathbf{1}\\
{\rm P_s}^T\cdot \mathbf{1}=0\\
\sum_{h=1}^{t}{\bar{\rm P}}_b^h+{{\rm E}}_b^0\le {\bar{\rm E}}_b\\
\sum_{h=1}^{t}{\bar{\rm P}}_b^h+{{\rm E}}_b^0\ge 0,
\end{aligned}
\label{eq:MILP1}
\end{equation}
where the efficiency of the battery is assumed to be 1, and ${\rm E}_b^0$ denotes the initial energy stored in the battery. As ${\rm P_s}^T\cdot \mathbf{1}$ is equal to the net energy stored in the battery during the entire day, the second constraint in \eqref{eq:MILP1} keeps the final SOC of the battery on its initial value. 
\subsection{Static and Quasi-Dynamic MILP Strategies}
\subsubsection{Static MILP}
Static MILP considers loads to be either on or off, i.e. transient dynamics are neglected. In this case, $\mathbf{R}$ and $\mathbf{Q}$ in equation \fref{const_bigM} have to be equal which forces $\mathbf{U}$ to be in the set of $\{0,1\}$. Static MILP still allows optimizing the size of a given number of units, minimum up and down time constraints, and the battery size and battery schedule for a given solar power profile. 
\subsubsection{Quasi-Dynamic MILP}
The quasi-dynamic MILP strategy is more sophisticated than static MILP. Here the units are assumed to go from off state to half operational state and finally to the final state, and vice versa. So the loads will go from 0\% to 50\% and then to 100\% of unit size $x_i$. This case is more realistic since most large loads are ramped up over a period of time before they reach to steady state.  The quasi-dynamic strategy also yields all the variables in the static strategy; the only difference is that $\mathbf{U}$ is in the set $\{0,1/2,1\}$.

\section{Numerical Results}\label{sec:results}
\begin {figure}[ht]
\graphicspath{ Plots/ }
\includegraphics[width=1\columnwidth, height=.6\columnwidth]{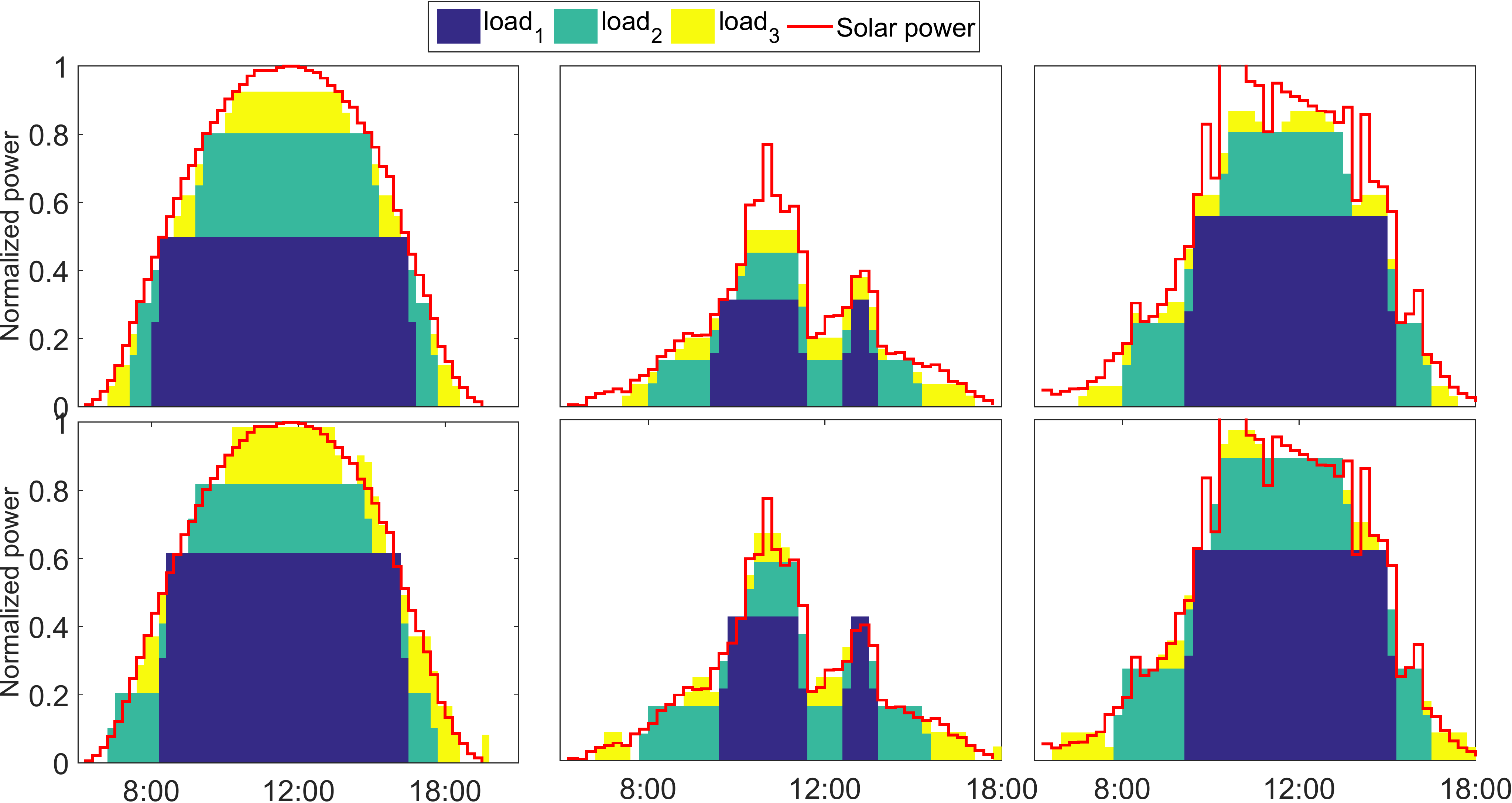}
\centering
\caption{Results of load sizing and scheduling without (top) and with (bottom) storage for the Quasi-Dynamic MILP strategy  for three different days and three units. The clear day (first from left) is referred to D1, the overcast day in the middle is D2, and the partially cloudy day (right) is D3. The minimum on and off times are 3 time steps (45 minutes) for all units. The initial storage of the battery is  \textbf{$\frac{\bar{P_b}}{2}$} and equal to [0.07, 0.13, 0.14]/2 for D1, D2, and D3 respectively}
%\jan also specify initial storage SOC. The x axis in between top and bottom in the center column do not appear to be teh same.
\label{figure:daysexamples}\vspace{-10pt}
\end {figure}

% \begin{figure*}[ht]\centering
% \includegraphics[width=2\columnwidth, height=.655\columnwidth,keepaspectratio]{drawing}
% \caption{Load scheduling for different solar forecast.}
% \label{Compare}\vspace{-15pt}
% \end{figure*}
In this section results from both static and quasi-dynamic MILP strategies are presented using using the CVX toolbox and Gurobi 6.50. Figure~\ref{figure:daysexamples} illustrates the results corresponding to three different representative solar power profiles of D1 clear day, D2 cloudy day and D3 partly cloudy day. The units are scheduled and stacked as area plots. The numeric values for units and battery are shown in Table~\ref{Semidynamicsresults}. Without storage (upper row in the figure) the total unit power cannot exceed the solar power. With storage (lower row) the total unit power can temporarily exceed solar power.

The optimal unit sizes clearly decrease on the overcast day with small solar resource as units of the size of $x_1$ in D1 would rarely be able to run. The addition of storage also effects the unit sizes. Storage will allow the units to ride through short-lived dips in solar power by utilizing battery energy which results in larger unit sizes and increased solar utilization. This is most apparent in the D3 case at around 1100 Local Standard Time (LST) where the green units ride through the temporary cloud cover on partial battery power and both blue and green units are enlarged for the storage case.  

In both Table~\ref{nodynamicsresults} and Table~\ref{Semidynamicsresults} different unit sizes and efficiency for the case $n=3$ are shown with different minimum up time and down time constraints. The left hand side column shows different minimum uptime and downtime constraints for each unit represented in terms of time steps in square brackets. At the right of the same column day types are shown as D1, D2, and D3. The middle column shows results for unit sizes and solar energy utilization or efficiency (Eff) without battery. Finally the last column displays the results with storage and the optimal battery size \textbf{$\bar{P}_b$} required to achieve Eff equal to one which corresponds to utilizing all solar power on a given day. All results including the battery sizes are normalized as per unit (pu) of the maximum solar power production on the clear day. Eff equal one can be achieved as the battery model is assumed ideal without losses from charging or discharging. 

\begin{table}[ht]
%\scriptsize
\centering
\caption{Static MILP strategy results ($\rm pu$) for three units and different on and off times (left column)
%\jan I am not sure about the objective with the different on and off times. I think the following would be more enlightening as to the effects of flexibility: [3 3 3] [3 3 3], [5 5 5] [5 5 5], [3 3 1] [3 3 1]
% The variations of minimum on versus off times for the same case do not seem to yield any patterns.
}
\label{nodynamicsresults}
\begin{tabular}{p{.4cm}p{.5cm}p{.3cm}|p{.3cm}p{.3cm}p{.3cm}p{.5cm}|p{.3cm}p{.3cm}p{.3cm}p{.3cm}}
                                                          \multirow{2}{*}{$\begin{bmatrix}m^+_1\\ m^+_2\\ m^+_3 \end{bmatrix} $}              &\multirow{2}{*}{$\begin{bmatrix}m^-_1\\ m^-_2\\ m^-_3 \end{bmatrix} $}                                                                       &    
                                                          \multirow{2}{*}{\begin{tabular}[c]{@{}l@{}}Day\\ type\end{tabular}} 
                                                          & \multicolumn{4}{c|}{Without Storage}                                                                  & \multicolumn{4}{c}{With Storage}                                                                         \\
%{$\begin{bmatrix}m^+_1\\ m^+_2\\ m^+_3 \end{bmatrix} $} 
& 
%{$\begin{bmatrix}m^-_1\\ m^-_2\\ m^-_3 \end{bmatrix} $} 
& %ype
%\multirow{3}{*}\begin{tabular}[c]{@{}l@{}}Day\\ type\end{tabular}  
&
$x_1$ & $x_2$ & $x_3$ & 
\begin{tabular}[c]{@{}l@{}}Eff\\ (\%)\end{tabular}  
& $x_1$ & $x_2$ & $x_3$ &
%begin{tabular}[c]{@{}l@{}}Battery\\ size\end{tabular}  
\textbf{$\bar{P}_b$} \\&    &      &      &      &                                                    &      &      &      &       
\\
%\hspace{.5cm}
%\hline 
\hhline{===|====|====}                             
&    &      &      &      &                                                    &      &      &      &        \\
    \multirow{3}{*}{$\begin{bmatrix}3\\ 3\\ 3 \end{bmatrix} $}                                   & \multirow{3}{*}{$\begin{bmatrix}3\\ 3\\ 3 \end{bmatrix} $}& D1& 0.50  & 0.30  & 0.12  & 0.90 & 0.54  & 0.33  & 0.12  & 0.11 \\
&&	D2	&	0.32	&	0.14	&	0.07	&	0.82	&	0.40	&	0.16	&	0.08	&	0.18	\\
&&	D3&	0.56	&	0.24	&	0.09	&	0.87&	0.62	&	0.27	&	0.11	&	0.13	\\		
&&&&&&&&&\\ \hline &&&&&&&&&\\\multirow{3}{*}{$\begin{bmatrix}3\\ 2\\ 1 \end{bmatrix} $}& \multirow{3}{*}{$\begin{bmatrix}3\\ 2\\ 1 \end{bmatrix} $}& D1& 0.55  & 0.25  & 0.12  & 0.92& 0.57  & 0.30  & 0.12  & 0.10\\
&&	D2&	0.29	&	0.15	&	0.11	&	0.84	&	0.38	&	0.18	&	0.12	&	0.12\\
&&	D3&	0.56	&	0.24	&	0.09	&	0.88	&	0.64	&	0.23	&	0.13	&	0.10\\
&&&&&&&&&\\ \hline &&&&&&&&&\\
\multirow{3}{*}{$\begin{bmatrix}7\\ 6\\ 5 \end{bmatrix} $}& \multirow{3}{*}{$\begin{bmatrix}3\\ 2\\ 1 \end{bmatrix} $}& D1& 0.60  & 0.30  & 0.08  & 0.89& 0.57  & 0.31  & 0.11  & 0.11\\
&&	D2&	0.25	&	0.20	&	0.11	&	0.81	&	0.33	&	0.24	&	0.10	&	0.18\\
&&	D3&	0.58	&	0.22	&	0.08	&	0.86	&	0.67	&	0.23	&	0.10	&	0.16\\
&&&&&&&&&\\ \hline &&&&&&&&&\\ 
\multirow{3}{*}{$\begin{bmatrix}3\\ 2\\ 1 \end{bmatrix} $}& \multirow{3}{*}{$\begin{bmatrix}7\\ 6\\ 5 \end{bmatrix} $}& D1& 0.50  & 0.30  & 0.16  & 0.91& 0.52  & 0.32  & 0.14  & 0.10\\
&&	D2&	0.25	&	0.20	&	0.11	&	0.81	&	0.55	&	0.22	&	0.11	&	0.19\\
&&	D3&	0.56	&	0.24	&	0.13	&	0.84	&	0.54	&	0.30	&	0.19	&	0.17
\end{tabular}
\end{table}

For three time steps (45~min) up and down time Table~\ref{nodynamicsresults} shows that clear and partly cloudy days (D1 and D3) were associated with larger units compared to the overcast day D2, as expected. Maybe unexpectedly, D3 had larger $x_1$ units compared to D1 but $x_{2,3}$ were smaller. D2 had smaller units to align with the smaller range of the solar resource and significant cloud variability. Efficiency was largest for the clear day as expected. 

Four different minimum uptime and downtime combinations were performed for 3 different daily patterns to illustrate the model sensitivity to different constraints. Larger minimum uptime and down time reduce the flexibility of scheduling and are expected to reduce efficiency and trigger larger unit sizes. When flexibility increases (from [3 3 3] to [3 2 1]) efficiency increases by about two percentage points and the required battery size for Eff$ = 1$ decreases by 10\% to 60\%. The results for different minimum on and off times for the same case are inconclusive and dependent on the time scales of solar resource fluctuations on different days.

Adding storage always increase Eff as our model solves for the co-optimization problem for both battery size and units size. All unit sizes increase when batteries are added especially on the more variables days D2 and D3. For all minimum uptime and downtime cases D2 requires the largest battery size to smooth out the variability even though the total solar production on D2 is far smaller than on D1 or D3. 
\begin{table}[ht]
\centering
\caption{Quasi-dynamic model results ($\rm pu$)
%\jan remove the last two rows here.
}
\label{Semidynamicsresults}
\begin{tabular}{p{.4cm}p{.5cm}p{.3cm}|p{.3cm}p{.3cm}p{.3cm}p{.5cm}|p{.3cm}p{.3cm}p{.3cm}p{.3cm}}
                                                                            \multirow{2}{*}{$\begin{bmatrix}m^+_1\\ m^+_2\\ m^+_3 \end{bmatrix} $}              &\multirow{2}{*}{$\begin{bmatrix}m^-_1\\ m^-_2\\ m^-_3 \end{bmatrix} $}                                                                       &      \multirow{2}{*}{\begin{tabular}[c]{@{}l@{}}Day\\ type\end{tabular}}& \multicolumn{4}{c|}{Without Storage}                                                                  & \multicolumn{4}{c}{With Storage}\\& & &$x_1$ & $x_2$ & $x_3$ & 
\begin{tabular}[c]{@{}l@{}}Eff\\ (\%)\end{tabular}  
& $x_1$ & $x_2$ & $x_3$ &
%begin{tabular}[c]{@{}l@{}}Battery\\ size\end{tabular}  
\textbf{$\bar{P}_b$}  \\&    &      &      &      &                                                    &      &      &      &       
\\
%\hspace{.5cm}
%\hline 
\hhline{===|====|====}                             
&    &      &      &      &                                                    &      &      &      &        
\\
\multirow{3}{*}{$\begin{bmatrix}3\\ 3\\ 3 \end{bmatrix} $}& \multirow{3}{*}{$\begin{bmatrix}3\\ 3\\ 3 \end{bmatrix} $}& D1 & 0.50 & 0.30 & 0.12 & 0.93& 0.54 & 0.32 & 0.12 & 0.07\\
&&	D2	&	0.32	&	0.14	&	0.07	&	0.85&	0.42	&	0.18	&	0.08	&	0.13\\
&&	D3	&	0.56	&	0.24	&	0.07	&	0.87&	0.71	&	0.22	&	0.11	&	0.14\\
&&&&&&&&&&\\
\hline &    &      &      &      &                                                    &      &      &      &        
\\
\multirow{3}{*}{$\begin{bmatrix}3\\ 2\\ 1 \end{bmatrix} $}                                   & \multirow{3}{*}{$\begin{bmatrix}3\\ 2\\ 1 \end{bmatrix} $}                    & D1 & 0.54 & 0.30 & 0.13 & 0.94                                               & 0.61 & 0.27 & 0.11 & 0.05                                                         \\
&&	D2	&	0.34	&	0.15	&	0.11	&	0.86 &	0.39	&	0.20	&	0.10	&	0.10\\
&&	D3	&	0.56	&	0.24	&	0.09	&	0.88&	0.57	&	0.29	&	0.14	&	0.10\\
                                    &                     &    &      &      &      &                                                    &      &      &      &                                                             \\ \hline                 &                     &    &      &      &      &                                                    &      &      &      &                                                             \\ 
                                                                                      \multirow{3}{*}{$\begin{bmatrix}7\\ 6\\ 5\end{bmatrix} $}                                   & \multirow{3}{*}{$\begin{bmatrix}3\\ 2\\ 1 \end{bmatrix} $}                         & D1 & 0.67 & 0.18 & 0.13 & 0.92                                               & 0.51 & 0.36 & 0.13 & 0.10                                                         \\
&&	D2	&	0.40	&	0.15	&	0.11	&	0.83&	0.34	&	0.25	&	0.12	&	0.18\\
&&	D3	&	0.56	&	0.24	&	0.06	&	0.86&	0.68	&	0.22	&	0.08	&	0.14\\
                                   &                      &    &      &      &      &                                                    &      &      &      &                                                             \\
                                    \hline &    &      &      &      &                                                    &      &      &      &        \\
                                                                                     \multirow{3}{*}{$\begin{bmatrix}3\\ 2\\ 1 \end{bmatrix} $}                                   & \multirow{3}{*}{$\begin{bmatrix}7\\ 6\\ 5 \end{bmatrix} $}                         & D1 & 0.55 & 0.37 & 0.25 & 0.91                                               & 0.63 & 0.35 & 0.18 & 0.10                                                         \\
&&	D2	&	0.41	&	0.16	&	0.09	&	0.82&	0.35	&	0.24	&	0.13	&	0.17\\
&&	D3	&	0.59	&	0.36	&	0.24	&	0.84&	0.54	&	0.37	&	0.20	&	0.17
\end{tabular}
\end{table}

Applying different minimum up and down times for different units are shown in the second row of results; the largest unit still has a 3 time step requirement, but the smallest unit does not have any constraint. The unit sizes for D3 remain largely unaffected, while for D1 $x_1$ is in fact larger and $x_2$ correspondingly smaller and vica versa for D2. The main noticeable difference here is the battery size of D2 which resulted in smaller size and due to the smaller minimum up/downtime constraints for the smaller units $x_{2}$ and $x_{3}$.

The last 2 cases are distinguishing both uptime downtime and units as shown in row 3 of the results, where the largest units $x_3$ has to be on for 7 times step (105~ min) but downtime is only 45~min, and so on. These cases are based on the load application assuming the units need to be on for double or more the time that they need to be off. Effects are mostly restricted to a redistribution of the unit sizes (e.g. for D2 $x_2$ increases while $x_1$ decreases) while the total unit sizes $x_1 + x_2 + x_3$ changes only slightly. Eff generally decreases with the reduced flexibility. In terms of storage, increasing the on time requirement did not effect much the D1 results, while the storage size increase substantially for D2 and D3. Lastly, flipping the previous cases assuming the load required longer off time, the effect is mainly in the unit sizes of D1 and D3. The Eff and \textbf{$\bar{P}_b$} remains similar as expected because the same solar profile results in similar battery scheduling.

Table \ref{Semidynamicsresults} shows results for the Quasi-Dynamic MILP strategy. The main change from the Static MILP results is increasing Eff in almost all cases since the half-on units effectively add additional discretizations beyond the $2^3$ unit size combinations for the three units. Therefore the loads can better track the solar curve. For the same reason \textbf{$\bar{P}_b$} became smaller. In general, unit sizes stayed the same or increased and in some cases like the last case in D1 ([3, 2, 1] and [7, 6, 5]) no small units resulted. For the clear day D1 usually $x_3$ is small as the middle of the day favors large base-load and the two smaller units primarily capture the shoulders of the day. However, long downtime constraints appear to favor larger units probably because the units are unable to capture the evening shoulder after the downtime requirement. 
\begin {figure}[ht]
\graphicspath{ Plots/ }
\includegraphics[width=1\columnwidth]{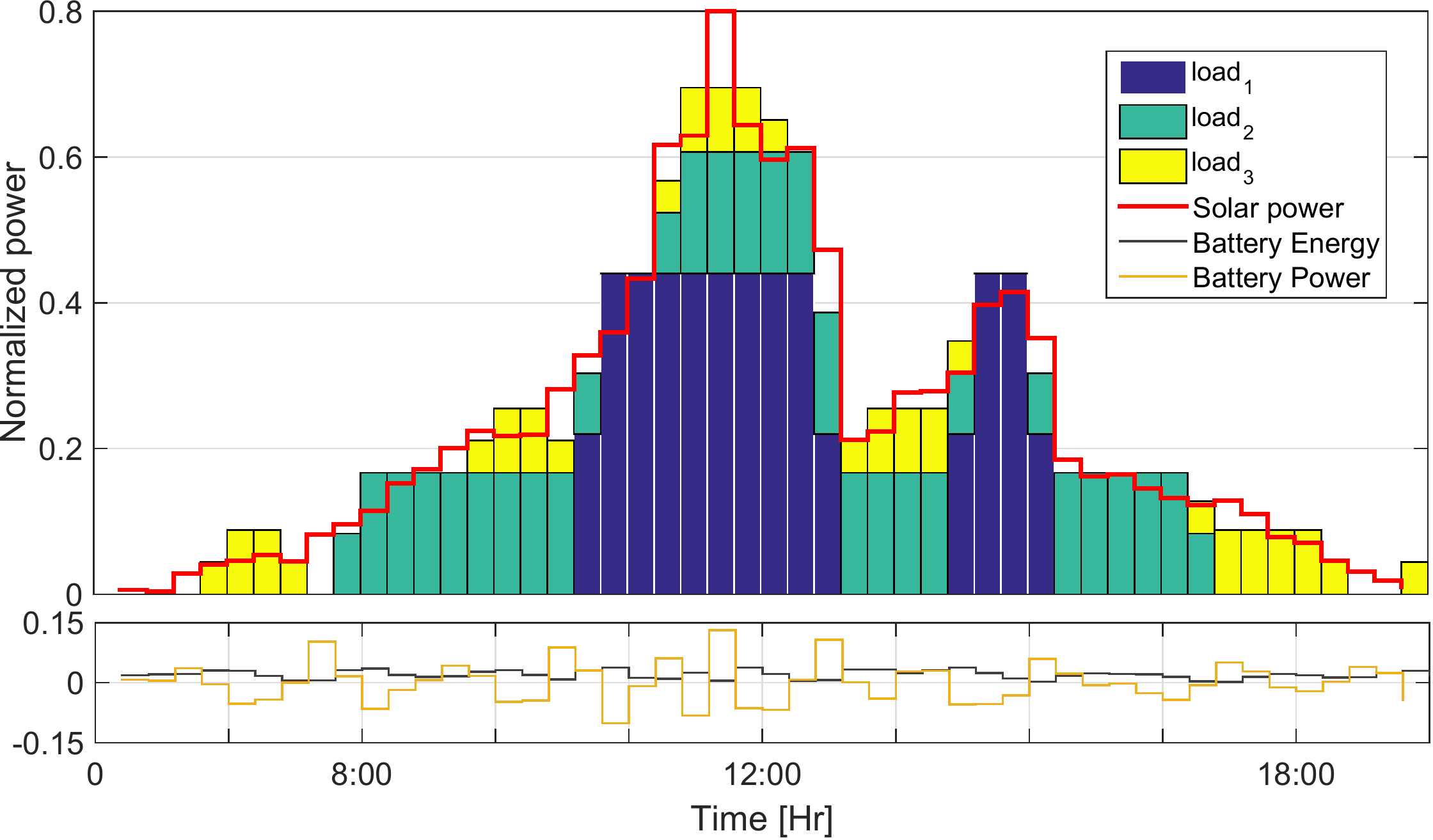}
\centering
\caption{Overcast day example of load sizing and scheduling with storage for quasi-dynamic MILP (top). Same as the bottom center graph in Fig. \ref{figure:daysexamples}, but the power and the energy of the battery are shown as well (bottom).
%\jan I think the legend for battery power and energy are flipped. In any case I am not able to see the black line. Remove the grid lines as the graph is already very busy. x axis seems to be time steps rather than hours.
}
\label{figure:Batterysizing}\vspace{-8pt}
\end {figure}

An overcast day example of load sizing and scheduling with storage is given in Figure~\ref{figure:Batterysizing}. The initial energy stored in the battery ${\rm E}_b^0$ is half of the battery size which is 0.13 pu and an additional constraint was added to the end charge state of the battery ${\rm E}_b^t$ = ${\rm E}_b^0$ or (${\rm P_s}\le {\bar{\rm P}}_b\cdot \mathbf{1}$). For the first few time steps after sunrise the battery charged allowing the smallest unit $x_3$ to be on for one hour. Then $x_3$ was shut down allowing the battery to charge up allowing $x_2$ to be 
\begin{table}[ht]
\centering
\caption{Different number of units ($n$) comparison for static MILP.}
\label{table:differentn}
\begin{tabular}{p{.4cm}p{.4cm}p{.4cm}p{.4cm}p{.4cm}p{.4cm}p{.5cm}p{.4cm}p{.4cm}p{.4cm}p{.4cm}}%{llllllllll}
\multicolumn{10}{c}{D1}                                              \\
$n$    & 2    & 3    & 4    & 5    &      & 2    & 3    & 4    & 5     \\
\textbf{$\bar{P}_b$}& 0    & 0    & 0    & 0    &      & 0.28  & 0.11  & 0.07  & 0.05  \\
$x_1$ & 0.67 & 0.50 & 0.48 & 0.49 &      & 0.67 & 0.54 & 0.49 & 0.50  \\
$x_2$ & 0.25 & 0.30 & 0.32 & 0.18 &      & 0.28 & 0.33 & 0.33 & 0.32  \\
$x_3$ &      & 0.12 & 0.12 & 0.14 &      &      & 0.12 & 0.13 & 0.20  \\
$x_4$ &      &      & 0.06 & 0.12 &      &      &      & 0.05 & 0.08  \\
$x_5$ &      &      &      & 0.05 &      &      &      &      & 0.07  \\
Eff     & 0.82 & 0.90 & 0.93 & 0.94 &      & 1 &    1  & 1 & 1  \\
\hhline{==========}     
\multicolumn{10}{c}{D2}                                              \\
$n$	&	2	&	3	&	4	&	5	&	&	2	&	3 &	4	&	5	\\
\textbf{$\bar{P}_b$}	&	0	&	0	&	0	&	0	&	&	0.40	&	0.18	&	0.09	&	0.05	\\
$x_1$	&	0.31	&	0.41	&	0.31	&	0.31	&	&	0.45	&	0.40	&	0.36	&	0.23	\\
$x_2$	&	0.14	&	0.14	&	0.15	&	0.26	&	&	0.15	&	0.16	&	0.22	&	0.21	\\
$x_3$	&	&	0.07	&	0.08	&	0.15	&	&	&	0.08	&	0.15	&	0.16	\\		
$x_4$	&	&	&	0.05	&	0.08	&	&	&	&	0.06	&	0.12	\\				
$x_5$	&	&	&	&	0.05	&	&	&	&	&	0.05	\\						
Eff	&	0.67	&	0.82	&	0.89	&	0.93	&	&	1	&	1	&	1&	1	\\
\hhline{==========}	
\multicolumn{10}{c}{D3}	\\																	
$n$	&	2	&	3 &	4	&	5	&	&	2	&	3	&	4	&	5	\\
\textbf{$\bar{P}_b$}	&	0	&	0	&	0	&	0	&	&	0.36	&	0.13	&	0.08	&	0.07	\\
$x_1$	&	0.56	&	0.56	&	0.56	&	0.40	&	&	0.69	&	0.62	&	0.58	&	0.51	\\
$x_2$	&	0.24	&	0.24	&	0.24	&	0.36	&	&	0.24	&	0.27	&	0.28	&	0.36	\\
$x_3$	&	&	0.09	&	0.24	&	0.18	&	&	&	0.11	&	0.12	&	0.18	\\		
$x_4$	&	&	&	0.09	&	0.08	&	&	&	&	0.06	&	0.08	\\				
$x_5$	&	&	&	0.04	&	0.06	&	&	&	&	&	0.06	\\					
	Eff	&	0.80	&	0.87	&	0.91	&	0.93	&	&	1	&	1	&	1	&	1
\end{tabular}
\end{table}
turned on 30~min before the solar power was sufficient. Following the charging and discharging battery curves shows that the battery was never charging or discharging for more than four consecutive time steps. The battery shifts small amounts of energy to allow units to come on earlier during increasing solar power production and turn off later during decreasing solar power production. 

Table~\ref{table:differentn} summarizes the results for different $n$ from two units up to five units to illustrate the ability of the model to simulate a varying number of units as well the size of each unit. The left side shows the sizing without storage based on optimal Eff. Eff increased from 82\% for two units up to 94\% for five units in the clear case and from 67\% to 93\% in the overcast case. On the right side storage is considered. For the clear day the system reaches optimal efficiency by using 28\% pu battery size for two units and reduces up to 5\% pu in the case of 5 units. For the overcast storage sizes decrease from 40\% to 5\%. 

The optimization code was performed in a 3.4 GHz Intel Core i7 processor with 32 GB of RAM. The computational time varies based on the day and the number of units, for the case of $n=2$ and an overcast day the computational time for the static method takes 2.5 seconds while quasi-dynamic method takes 10 seconds the computational time. 
%For a yearly results, an option of running all day and then by statistical approaches determine with is the best sizing and number of unis. %Another option is clustering how many clear or cloudy day there in the a selected location and weight the results of each day.
 
\section{Conclusions}\label{sec:CCLS}
Considering the intermittency of renewable energy systems, a sizing and scheduling model is proposed for a finite number of static or quasi-dynamic electric loads. The model objective is to maximize system efficiency, which is also defined as solar utilization, with and without storage. For the application of optimal load size selection, the power production of a solar PV is assumed to be consumed by a finite number of discrete loads in an off-grid system using mixed-integer linear programming with constraints, such as, battery charge and discharge limits, and minimum uptime and downtime for each unit. The method was applied to three characteristic daily solar profiles. Different minimum up and down time constraints are also investigated. 

By means of a case study, three different days results indicate a system efficiency increased from 82\% for two units up to 94\% for five units for the clear day and from 67\% to 93\% for the overcast day. Including battery storage for the clear day, the system requires a 28\% pu battery size to reach 100\% efficiency for two units, but the battery size reduces to 5\% pu for five units.

The results obtained are specific to the location and days presented differ by location and may even vary year-to-year due to spatio-temporal patterns in the solar resources and clouds coverage. The methodology proposed in this paper allows computationally efficient solutions even when several years of solar resource data are available and yield the optimal sizing for the given data. 

For practical applications, the economics also need to be considered as smaller units typically cost more per kW and an optimization based on cost would therefore yield larger and fewer units. Within our framework, it is possible to assign a cost function to the number of units and to the efficiency to allow satisfying needs of practitioners. Similarly, the competition between reduced battery size and larger unit capital cost for more units could be considered in such an economic optimization.

% The model can also be applied to derive the optimum number of units. While for the present cost function the optimal number of units is infinity, overall economic costs tend to increase when the same load is split between many smaller units. If this behavior was described through an economic cost function than the model could solve for the optimal number of units. 

\bibliographystyle{IEEEtran}
\bibliography{mylib}
\end{document}